# THE FOUR VERTEX THEOREM AND ITS CONVERSE

*Dennis DeTurck, Herman Gluck, Daniel Pomerleano, David Shea Vick*

**Dedicated to the memory of Björn Dahlberg**

## Abstract

The *Four Vertex Theorem*, one of the earliest results in global differential geometry, says that a simple closed curve in the plane, other than a circle, must have at least four "vertices", that is, at least four points where the curvature has a local maximum or local minimum. In 1909 Syamadas Mukhopadhyaya proved this for strictly convex curves in the plane, and in 1912 Adolf Kneser proved it for all simple closed curves in the plane, not just the strictly convex ones.

The *Converse to the Four Vertex Theorem* says that any continuous real-valued function on the circle which has at least two local maxima and two local minima is the curvature function of a simple closed curve in the plane. In 1971 Herman Gluck proved this for strictly positive preassigned curvature, and in 1997 Björn Dahlberg proved the full converse, without the restriction that the curvature be strictly positive. Publication was delayed by Dahlberg's untimely death in January 1998, but his paper was edited afterwards by Vilhelm Adolfsson and Peter Kumlin, and finally appeared in 2005.

The work of Dahlberg completes the almost hundred-year-long thread of ideas begun by Mukhopadhyaya, and we take this opportunity to provide a self-contained exposition.

August, 2006

# I. Why is the Four Vertex Theorem true?

## 1. A simple construction.

A counter-example would be a simple closed curve in the plane whose curvature is nonconstant, has one minimum and one maximum, and is weakly monotonic on the two arcs between them. We will try to build such a curve from a few arcs of circles, fail, and see why. Consider the figure below.

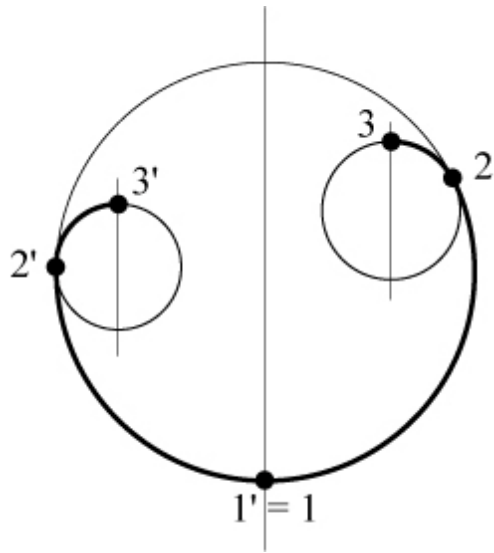

Think of the circles above as train tracks. Our curve will be traced out by a train which travels along these circular tracks, switching from one to another at the labeled junction points.

Let the train start at point 1 on the bottom of the largest circle, giving the minimum curvature there. It then moves off to the right along this large circle, comes to the point 2, and switches there to a smaller circle, increasing its curvature.

The important thing to note is that the vertical diameter of this smaller circle is displaced to the right of the vertical diameter of the original circle.

To keep things brief, we let the train stay on this second circle until it comes to the top at the point 3, which is to the right of the vertical diameter of the original circle.

Then we go back to the beginning, start the train once again at the bottom of the original circle, relabel this point 1', and now let the train move off to the left. It switches to a smaller circle at point 2', and comes to the top of this circle at 3', which is now to the left of the vertical diameter of the original circle.

If the curved path of the train is to form a convex simple closed curve, then the top points 3 and 3' should coincide...but they don't. So the curve does not exist.

No matter how many circles we use, we get the same contradiction: the curve can not close up while staying simple, and therefore does not exist.

But if we permit a self-intersection, then it's easy to get just one maximum and one minimum for the curvature.

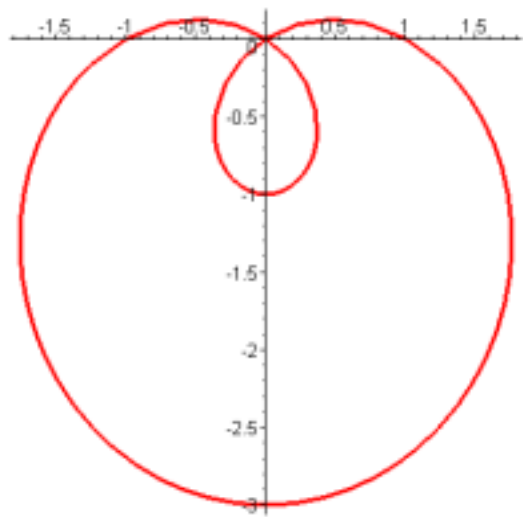

The equation of the above curve in polar coordinates is $r = -1 - 2 \sin \theta$. The curvature takes its minimum value of $5/9$ at the bottom of the big loop, its maximum value of $3$ at the bottom of the small loop, and is monotonic in between.

# II. Osserman's proof of the Four Vertex Theorem.

## 2. Consider the circumscribed circle.

In 1985, Robert Osserman gave a simple proof of the Four Vertex Theorem in which all cases, strictly positive or mixed positive and negative curvature, are treated on an equal footing. At the beginning of his paper, Osserman says,

> *"The essence of the proof may be distilled in a single phrase:
> consider the circumscribed circle."*

Let $E$ be a fixed nonempty compact set in the plane. Among all circles $C$ which enclose $E$, there is a unique smallest one, the ***circumscribed circle***.

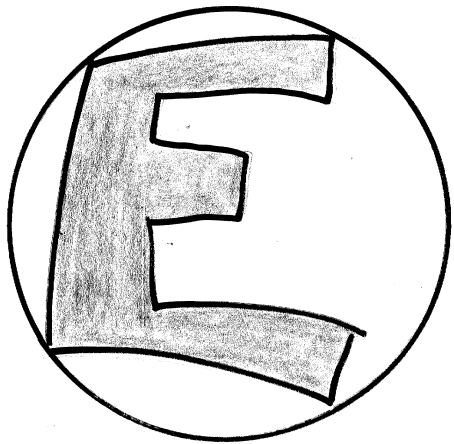

What can $C \cap E$ look like ?

(1) $C$ must meet $E$, otherwise $C$ could be made smaller and still enclose $E$.

(2) $C \cap E$ can not lie in an open semi-circle. Otherwise we could shift $C$ a bit, without changing its size, so that the new $C$ still encloses $E$ but does not meet it. But this would contradict (1).

(3) Thus $C \cap E$ has at least two points, and if only these, then they must be antipodal. This happens when $E$ is an ellipse.

## 3. An overview of Osserman's proof.

**OSSERMAN'S THEOREM.** *Let $\alpha$ be a simple closed curve of class $C^2$ in the plane, and C the circumscribed circle. If $\alpha \cap C$ has at least n components, then $\alpha$ has at least 2n vertices.*

Here are the major steps in the proof.

(1) Let R be the radius of the circumscribed circle C, and $K = 1/R$ its curvature. If $\alpha \cap C$ has at least n components, then $\alpha$ has at least n vertices where the curvature $\kappa$ satisfies $\kappa \geq K$, and at least n more vertices where $\kappa < K$.

*Caution.* We do not claim that the six points shown in the figure below are vertices. But if you find six consecutive points on the curve $\alpha$ where the curvature goes up and down as shown, then you can replace them by six vertices of $\alpha$, three of them local maxima and three local minima. We leave it to the reader to check this simple assertion.

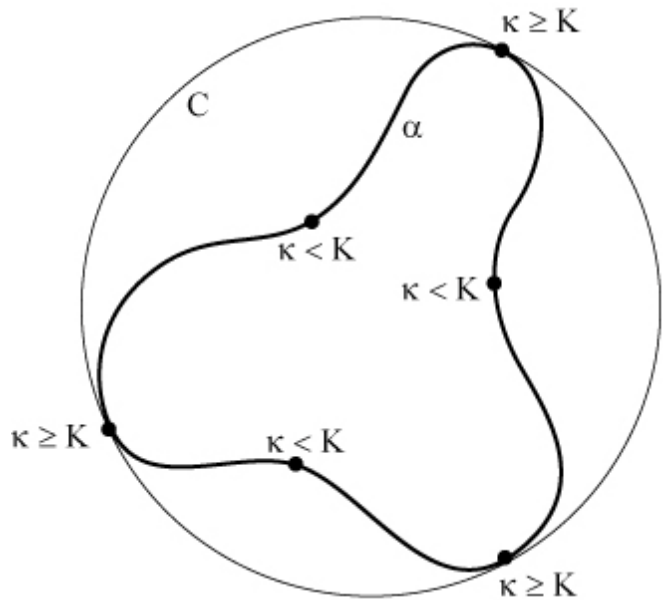

It is clear from this picture how we get the n points where $\kappa \geq K$.

To get the ones inbetween, where $\kappa < K$, we proceed as follows.

Let $P_1$, $P_2$, ..., $P_n$ be $n$ points in counterclockwise order on the circle $C$, one from each component of $C \cap \alpha$. Since $C \cap \alpha$ cannot lie in an open semi-circle, we can choose these points so that the arc of $C$ from each point to the next is no larger than a semi-circle.

Now focus on the arc $C_1$ of $C$ from $P_1$ to $P_2$, and on the corresponding arc $\alpha_1$ of the curve $\alpha$. We have arranged for the convenience of description that the line connecting $P_1$ and $P_2$ is vertical.

The arc $\alpha_1$ of $\alpha$ is tangent to $C$ at $P_1$ and $P_2$, and does not lie entirely along $C$. It therefore contains a point $Q$ which is not on $C_1$, yet lies to the right of the vertical line connecting $P_1$ and $P_2$. Consider the circle $C'$ through $P_1$, $Q$ and $P_2$.

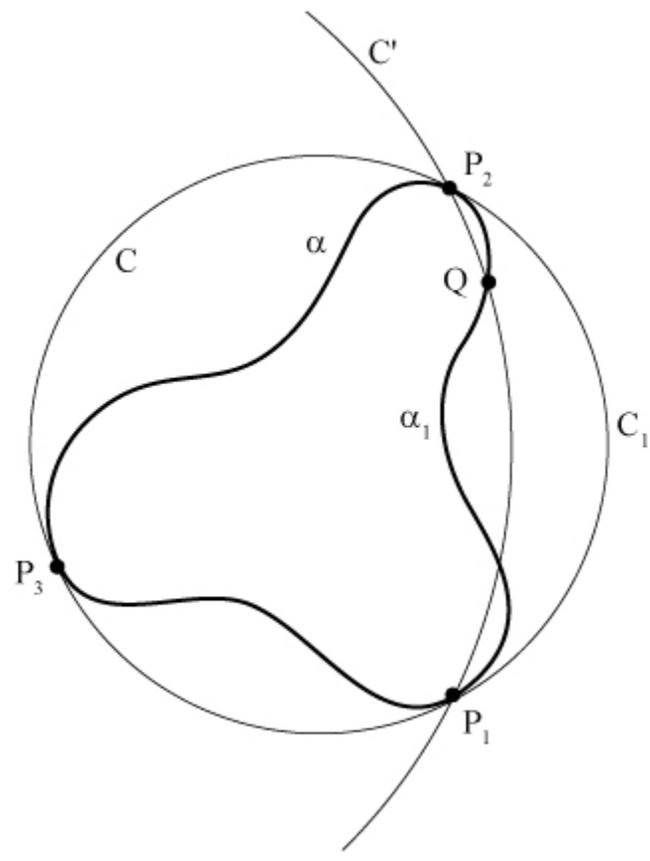

Since $Q$ is in the interior of $C$, and since $C_1$ is no larger than a semicircle, it follows that the circle $C'$ is larger than $C$. Hence its curvature $K' < K$.

Now gradually translate the circle C' to the left, until it last touches $\alpha_1$ at the point $Q_1$.

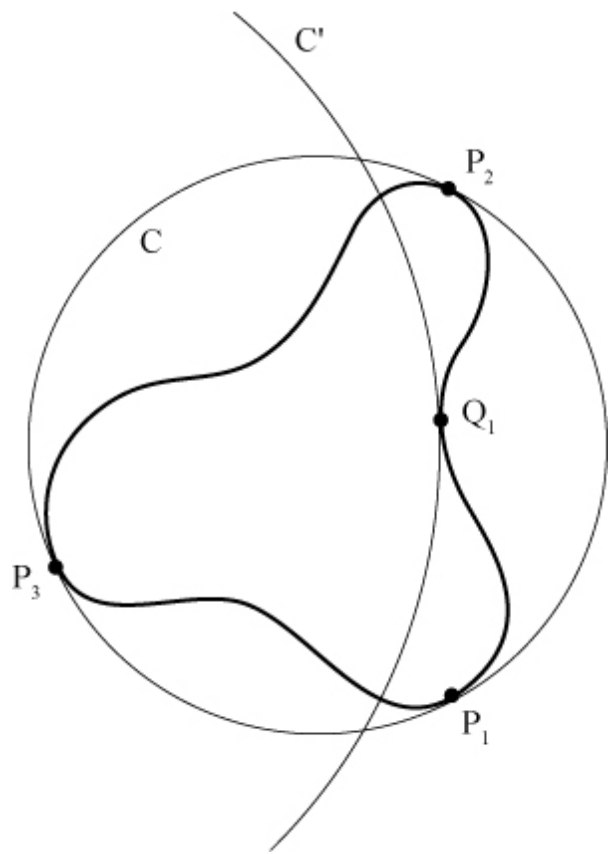

At $Q_1$, the curvature $\kappa$ of $\alpha$ must satisfy $\kappa(Q_1) \leq K' < K$, which is exactly what we wanted to show.

(2) **Bonus clause.** The components of $\alpha \cap C$ are either single points or closed arcs. For each component which is an arc rather than a single point, we can get two extra vertices, over and above the $2n$ vertices promised by Osserman's theorem.

The picture of $\alpha$ below is roughly the same as before, but we have distorted it so that an entire arc of $\alpha$ about the point $P_2$ coincides with an arc of the circle $C$. The curvature of $\alpha$ at $P_2$ is now exactly $\kappa$. But in addition, there must be points of $\alpha$ to either side of this arc, arbitrarily close to it, where the curvature $\kappa$ of $\alpha$ is $> K$. Two such points are shown in the figure, and marked $R_2$ and $R_2'$. As a result, the point $P_2$, where the curvature of $\alpha$ used to be larger than the curvature at the two neighboring marked points, has been demoted, and the curvature of $\alpha$ at $P_2$ is now less than that at the two neighboring marked points. In this way, we have gained two extra vertices for $\alpha$.

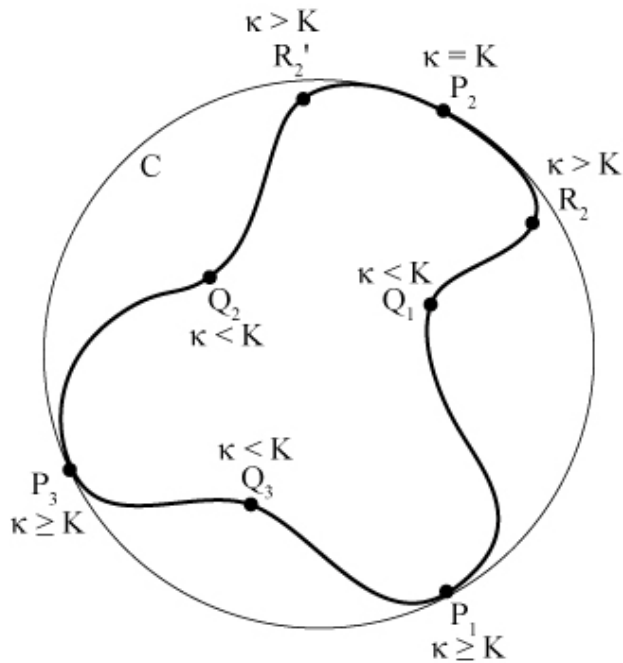

(3) It can happen that $\alpha \cap C$ has only one component, which must then be a closed arc on $C$ greater than or equal to a semi-circle. Osserman's theorem promises that $\alpha$ has at least two vertices, and the bonus clause promises two more.

(4) Thus Osserman's theorem, together with the bonus clause, yields the Four Vertex Theorem.

# III. The converse in the case of strictly positive curvature.

**FULL CONVERSE TO THE FOUR VERTEX THEOREM.** *Let $\kappa: S^1 \to R$ be a continuous function which is either a nonzero constant or else has at least two local maxima and two local minima. Then there is an embedding $\alpha: S^1 \to R^2$ whose curvature at the point $\alpha(t)$ is $\kappa(t)$ for all $t \in S^1$.*

## 4. Basic idea of the proof for strictly positive curvature.

We use a winding number argument in the group of diffeomorphisms of the circle, as follows.

Let $\kappa: S^1 \to R$ be any continuous strictly positive curvature function, and for the moment think of the parameter along $S^1$ as the angle of inclination $\theta$ of the desired curve. Since this curve will usually not close up, it is better to cut $S^1$ open into the interval $[0, 2\pi]$.

There is a unique map $\alpha: [0, 2\pi] \to R^2$ which begins at the origin and has unit tangent vector $(\cos \theta, \sin \theta)$ and curvature $\kappa(\theta)$ at the point $\alpha(\theta)$. In fact, if s is arc length along this curve, then the equations

$$d\alpha/ds = (\cos \theta, \sin \theta) \quad \text{and} \quad d\theta/ds = \kappa(\theta)$$

quickly lead to the explicit formula

$$\alpha(\theta) = \int_0^\theta (\cos \theta, \sin \theta) \, d\theta / \kappa(\theta) .$$

The *error vector* $E = \alpha(2\pi) - \alpha(0)$ measures the failure of our curve to close up.

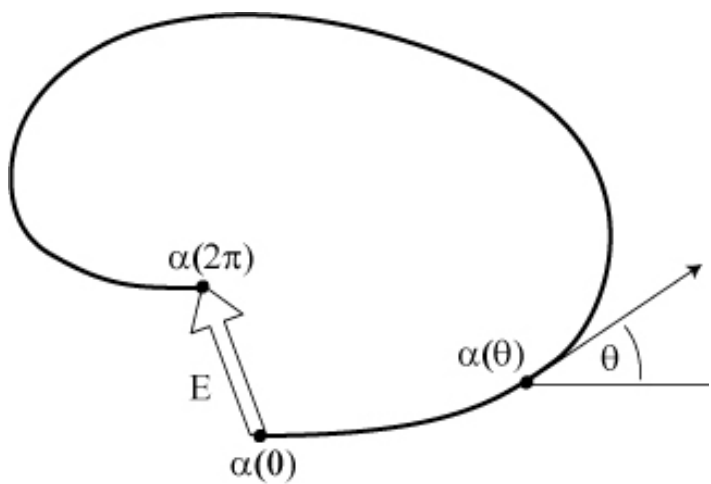

If the curvature function $\kappa$ has at least two local maxima and two local minima, we will show how to find a loop of diffeomorphisms $h$ of the circle so that when we construct, as above, curves $\alpha_h: [0, 2\pi] \to R^2$ whose curvature at the point $\alpha_h(\theta)$ is $\kappa \circ h(\theta)$, the corresponding error vectors $E(h)$ will wind once around the origin. Furthermore, we will do this so that the loop is contractible in the group $\text{Diff}(S^1)$ of diffeomorphisms of the circle, and conclude that for some $h$ "inside" this loop, the curve $\alpha_h$ will have error vector $E(h) = 0$, and hence close up to form a smooth simple closed curve.

Its curvature at the point $\alpha_h(\theta)$ is $\kappa \circ h(\theta)$, so if we write $\alpha_h(\theta) = \alpha_h \circ h^{-1} \circ h(\theta)$, and let $h(\theta) = t$, then its curvature at the point $\alpha_h \circ h^{-1}(t)$ is $\kappa(t)$. Therefore $\alpha = \alpha_h \circ h^{-1}$ is a reparametrization of the same curve, whose curvature at the point $\alpha(t)$ is $\kappa(t)$.

## 5. How do we find such a loop of diffeomorphisms?

The first step is to replace $\kappa$ by a simpler curvature function.

Using the hypothesis that $\kappa$ has at least two local maxima and two local minima, we find positive numbers $0 < a < b$ such that $\kappa$ takes the values $a, b, a, b$ at four points in order around the circle. We then find a preliminary diffeomorphism $h_1$ of the circle so that the function $\kappa \circ h_1$ is $\varepsilon$-*close in measure* to a step function $\kappa_0$ with values $a, b, a, b$ on four successive arcs of length $\pi/2$ each, in the sense that $\kappa \circ h_1$ is within $\varepsilon$ of $\kappa_0$ on almost all of $S^1$, except on a set of measure $< \varepsilon$.

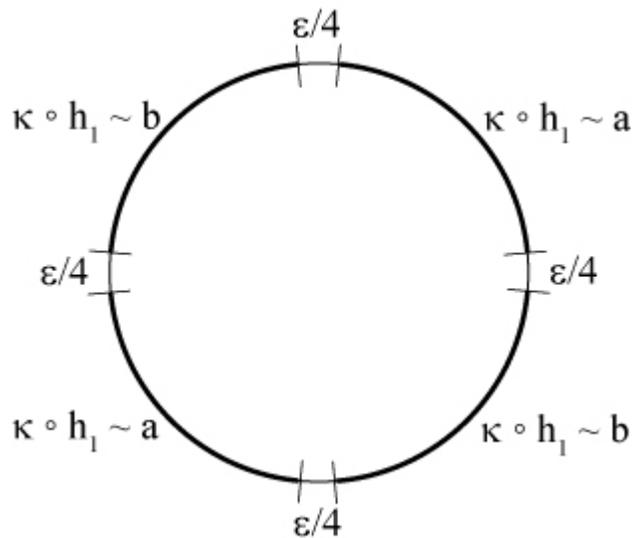

If the curvature functions $\kappa \circ h_1$ and $\kappa_0$ are $\varepsilon$-close in measure, then the corresponding curves are easily seen to be $C^1$-close. With this in mind, we will focus attention on the curvature step function $\kappa_0$ and ignore the original curvature function $\kappa$ until almost the end of the argument.

Given this shift of attention, we next find a contractible loop of diffeomorphisms $h$ of the circle so that the error vectors $E(h)$ for the curvature step functions $\kappa_0 \circ h$ wind once around the origin, as in the picture below.

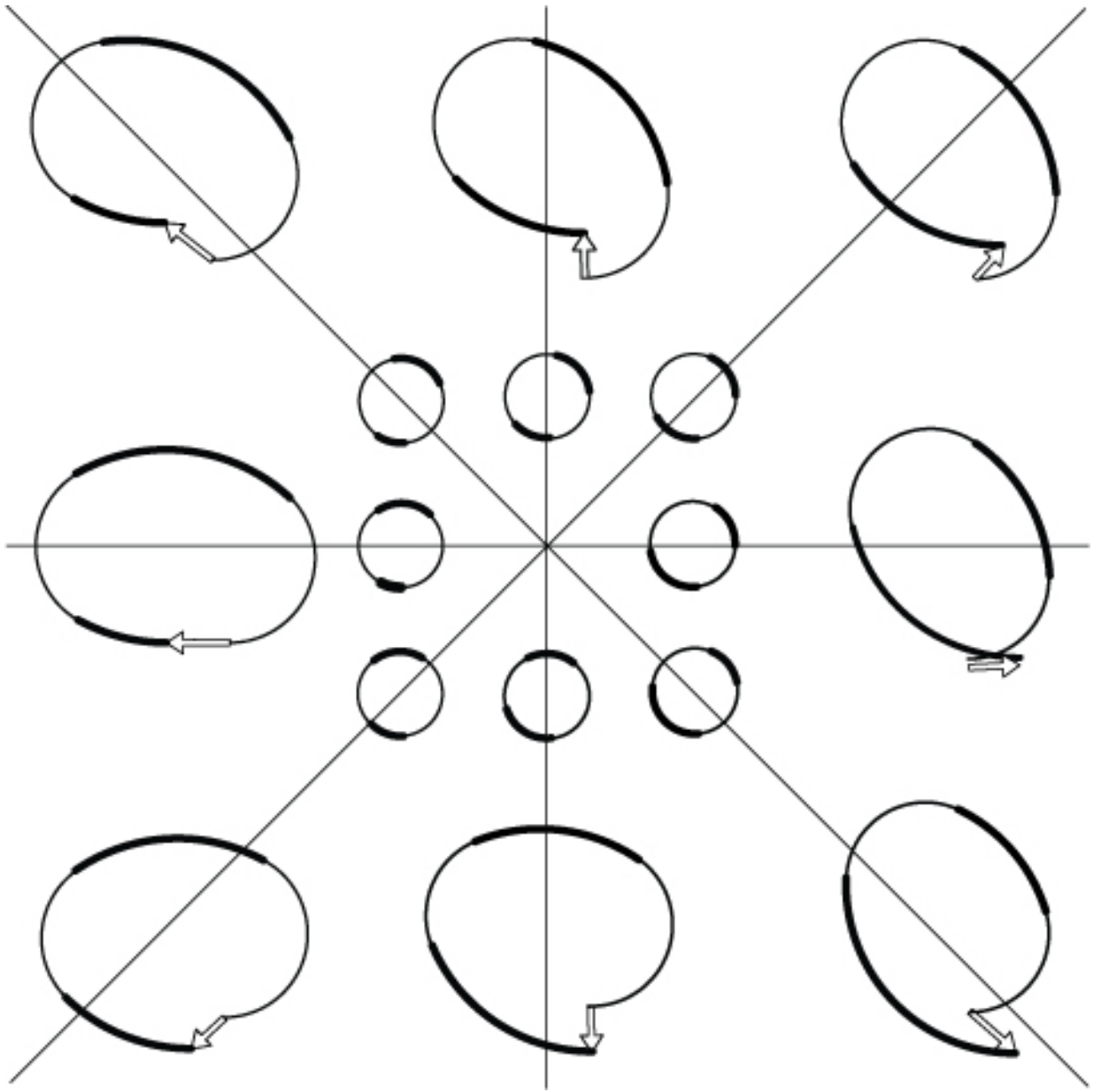

**A curve tries unsuccessfully to close up**

# 6. How do we read this picture?

The picture has eight portions, located at the eight points of the compass. Each contains a small circle which shows a curvature step function $\kappa_0 \circ h$, and a larger curve which realizes it, with common parameter the angle of inclination $\theta$. This larger curve is built from four circular arcs, two cut from a circle of curvature $a$ and two from a circle of curvature $b$. Heavier markings indicate the larger circle.

The curves fail to close up, and we show their error vectors. The figures are arranged so that the one in the east has an error vector pointing east, the one in the northeast has an error vector pointing northeast, and so on.

The eight diffeomorphisms $h$ are samples from a loop of diffeomorphisms of the circle, and the corresponding error vectors $E(h)$ wind once around the origin. The loop is contractible in $\text{Diff}(S^1)$ because each diffeomorphism leaves the bottom point of the circle fixed. The precise description of such a family of diffeomorphisms of the circle is given in Gluck (1971), but we think there is enough information in the above picture for the reader to work out the details. This picture can be seen in animated form at our website, *www.math.upenn.edu/~deturck/fourvertex*.

Since the loop of diffeomorphisms is contractible in $\text{Diff}(S^1)$, it follows that for some other diffeomorphism $h$ of the circle, the corresponding curve with curvature $\kappa_0 \circ h$ has error vector zero, and therefore closes up. This is hardly news: if $h$ is the identity, then the curve with curvature $\kappa_0$ closes up, as shown below.

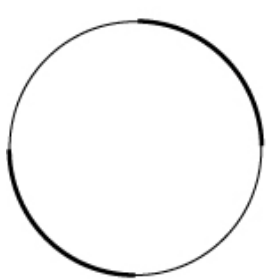 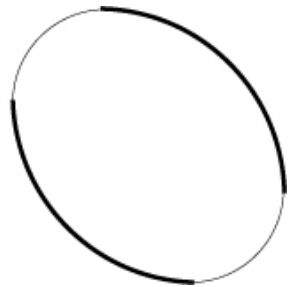

**Preassign $\kappa_0$**       **Get this "bicircle"**

But the point is that the above argument is *robust*, and hence applies equally well to the curvature function $\kappa \circ h_1$ which is $\varepsilon$-close in measure to $\kappa_0$, since the corresponding curves are $C^1$-close. Hence there must also be a diffeomorphism $h$ of the circle so that the curve with curvature $\kappa \circ h_1 \circ h$ has error vector zero, and therefore closes up to form a strictly convex simple closed curve. Reparametrizing this curve as in section 4 shows that it realizes the preassigned curvature function $\kappa$, completing the proof of the converse to the Four Vertex Theorem for strictly positive curvature.

## 7. Why does this argument sound familiar?

Proving something that is obvious by using a robust argument reminds us of the proof that the polynomial $p_0(z) = a_n z^n$ has a zero by noting that it takes any circle in the complex plane to a curve which winds n times around the origin. After all, it is obvious that $p_0(z)$ has a zero at the origin.

But this argument is also robust, and applies equally well to the polynomial

$$p(z) = a_n z^n + a_{n-1} z^{n-1} + \ldots + a_1 z + a_0,$$

when z moves around a suitably large circle, giving us the usual topological proof of the Fundamental Theorem of Algebra.

## 8. Where does our winding number argument come from?

Here is a higher-dimensional result proved in Gluck (1972).

**GENERALIZED MINKOWSKI THEOREM.** *Let* $K: S^n \to R$, *for* $n \geq 2$, *be a continuous, strictly positive function. Then there exists an embedding* $\alpha: S^n \to R^{n+1}$ *onto a closed convex hypersurface whose Gaussian curvature at the point* $\alpha(p)$ *is* $K(p)$ *for all* $p \in S^n$.

The proof of this theorem is by a degree argument in the diffeomorphism group of the n-sphere $S^n$. The hypothesis that $n \geq 2$ is used to guarantee the existence of a diffeomorphism which arbitrarily permutes any finite number of points. This hypothesis fails for $n = 1$, but when we rewrite the argument in this case, we get the winding number proof of the converse to the Four Vertex Theorem for strictly positive curvature described above.

## 9. Building curves from arcs of circles.

Before we leave this section, we record the following observation of Dahlberg, and save it for later use.

**PROPOSITION 9.1.** *A plane curve, built with four arcs cut in alternation from two different size circles, and assembled so that the tangent line turns continuously through an angle of $2\pi$, will close up if and only if opposite arcs are equal in length.*

Evidence for this proposition is provided by the preceding figure of a bicircle, and the large figure before that, containing eight separate curves which fail to close up.

**Proof.** We parametrize such a curve by the angle of inclination $\theta$ of its tangent line, suppose that $0 < a < b$ are the curvatures of the two circles to be used, and record the curvature function $\kappa(\theta)$ in the figure below.

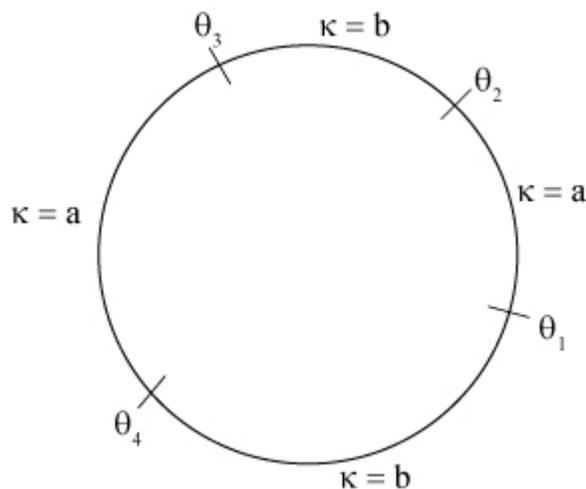

If $\alpha(\theta)$ is the resulting curve, then its error vector is given by

$$E = \alpha(2\pi) - \alpha(0) = \int_0^{2\pi} e^{i\theta}\, ds = \int_0^{2\pi} e^{i\theta}\, ds/d\theta\, d\theta = \int_0^{2\pi} e^{i\theta}/\kappa(\theta)\, d\theta.$$

We evaluate this explicitly as a sum of four integrals, and get

$$E = (1/ia - 1/ib)\left([\exp(i\theta_2) - \exp(i\theta_1)] + [\exp(i\theta_4) - \exp(i\theta_3)]\right).$$

We are left with the question of whether the two vectors shown below add up to zero.

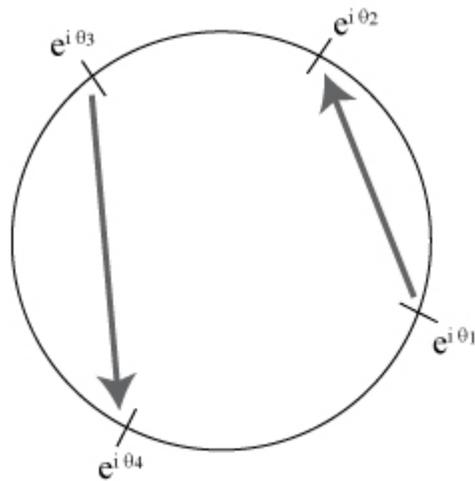

The picture makes the answer obvious and proves the proposition.

# IV. Dahlberg's proof of the Full Converse to the Four Vertex Theorem.

## 10. Dahlberg's key idea.

When asked to draw a simple closed curve in the plane with strictly positive curvature, and another one with mixed positive and negative curvature, a typical response might be...

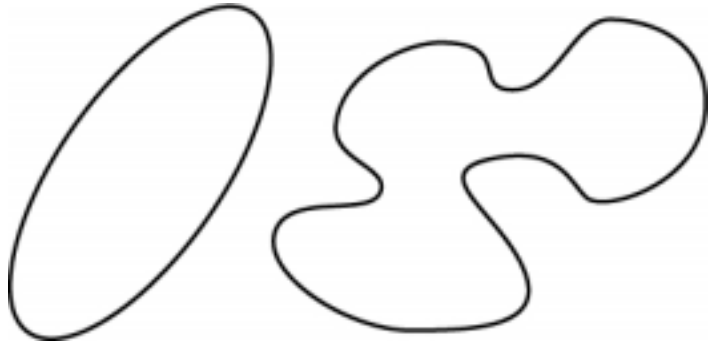

But for mixed positive and negative curvature, Dahlberg envisioned the following curve.

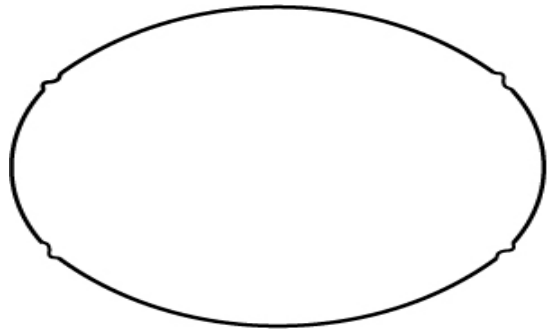

Its four major subarcs are almost circular. They are connected by four small wiggly arcs, each of which has an almost constant tangent direction, but largely varying curvature, including negative curvature. This curve "marginalizes" its negative curvature and emphasizes its positive curvature, and from a distance looks like a bicircle.

**Dahlberg's Key Idea.** *You can construct such a curve with any preassigned curvature which has at least two local maxima and two local minima. You can use the winding number argument to get it to close up smoothly, and you can also make it $C^1$ close to a fixed convex curve, which will imply that it is simple.*

## 11. Dahlberg's proof plan.

The overall plan is to follow the winding number line of argument from the strictly positive curvature case, preserving the prominent role played by two-valued curvature step functions, but mindful of the need for the following changes.

(1) A smooth simple closed curve with strictly positive curvature has a distinguished parametrization by the angle of inclination of its tangent line. This doesn't work in general, so another parametrization scheme is needed.

(2) A smooth closed curve with strictly positive curvature whose tangent line turns through an angle of $2\pi$ is automatically simple. This is false in general, and part of Dahlberg's Key Idea is to force the curve to be $C^1$ close to a fixed convex curve in order to make it simple.

To accomplish this, he exhibits a 2-cell **D** in $\text{Diff}(S^1)$, centered at the identity, and satisfying a certain transversality condition which guarantees that the winding number argument will work using arbitrarily small loops in **D** about its center.

## 12. Finding the right parametrization.

Choose, as common domain for all our curves, the unit circle $S^1$ with arc length s as parameter. Since most of the curves will fail to close up, cut the circle open at the point $(1,0)$ and think of the interval $[0, 2\pi]$ as the domain. All the curves will have length $2\pi$, and at the end will be scaled up or down to modify their curvature.

We make sure the total curvature is always $2\pi$ as follows. Given a preassigned curvature function $\kappa: S^1 \to R$ which is not identically zero, evaluate $\int_0^{2\pi} \kappa(s) \, ds$. If this is zero, precede $\kappa$ by a preliminary diffeomorphism of $S^1$ so as to make the integral nonzero, but still call the composition $\kappa$. Then rescale this new $\kappa$ by a constant c so that

$$\int_0^{2\pi} c \, \kappa(s) \, ds = 2\pi.$$

When we later modify this curvature function by another diffeomorphism $h: S^1 \to S^1$, we will rescale the new curvature function $\kappa \circ h$ by a constant $c_h$ so that the total curvature is again $2\pi$:

$$\int_0^{2\pi} c_h \, \kappa \circ h(s) \, ds = 2\pi.$$

We then build a curve $\alpha_h: [0, 2\pi] \to R^2$ parametrized by arc length which begins at the origin, $\alpha_h(0) = (0, 0)$, starts off in the direction of the positive x-axis, $\alpha_h'(0) = (1, 0)$, and whose curvature at the point $\alpha_h(s)$ is $c_h \, \kappa \circ h(s)$. The curve $\alpha_h$ is uniquely determined by these requirements.

The effect of these arrangements is that all our curves begin and end up pointing horizontally to the right. If any such curve closes up, it does so smoothly. Thus, just as in the positive curvature case, the emphasis is on getting the curve to close up.

At the end, using ideas expressed above, we will make sure the final curve is simple.

## 13. Configuration space.

Defining a configuration space of four ordered points on a circle will help us visualize the arguments to come, and aid us in carrying out the transversality arguments mentioned in the proof plan.

As in the positive curvature case, we use the hypothesis that $\kappa$ has at least two local maxima and two local minima to find positive numbers $0 < a < b$ so that $\kappa$ takes the values $a, b, a, b$ at four points in order around the circle. The only difference here is that we might first have to change the sign of $\kappa$ to accomplish this.

Then we find a preliminary diffeomorphism $h_1$ of the circle so that the function $\kappa \circ h_1$ is $\varepsilon$-close in measure to the step function $\kappa_0$ with values $a, b, a, b$ on the arcs $[0, \pi/2), [\pi/2, \pi), [\pi, 3\pi/2), [3\pi/2, 2\pi)$.

As before, we focus on the step function $\kappa_0$ and its compositions $\kappa_0 \circ h$ as $h$ ranges over $\text{Diff}(S^1)$. These are all step functions with the values $a, b, a, b$ on the four arcs determined by some four points $p_1, p_2, p_3, p_4$ in order on $S^1$. For each such step function, we follow Dahlberg's parametrization scheme from the previous section to construct a curve from circular arcs with curvatures proportional to $a, b, a, b$, scaled up or down to make the total curvature $2\pi$. Let $E(p_1, p_2, p_3, p_4)$ denote the error vector for this curve.

To think visually about all this, let CS denote the configuration space of ordered 4-tuples $(p_1, p_2, p_3, p_4)$ of distinct points on the unit circle $S^1$, arranged in counterclockwise order, as shown below.

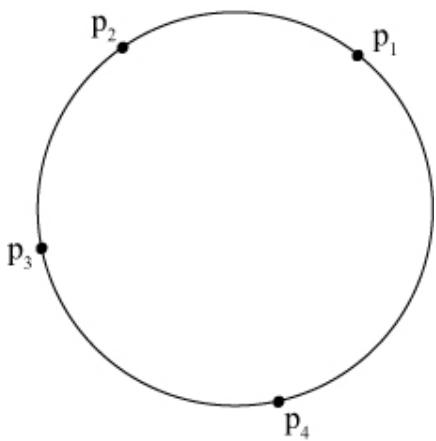

The configuration space CS is diffeomorphic to $S^1 \times R^3$.

The error vector $E(p_1, p_2, p_3, p_4)$ defines an *error map* $E: CS \to R^2$. The vanishing of the error vector is the signal that the curve closes up. We will call the set of such points the *core* of the configuration space $CS$, and denote it by $CS_0$.

We noted in Proposition 9.1 that a curve built with four arcs cut in alternation from two different size circles will close up if and only if opposite arcs are equal in length.

Thus a point $(p_1, p_2, p_3, p_4)$ of $CS$ lies in the core if and only if $p_1$ and $p_3$ are antipodal, and also $p_2$ and $p_4$ are antipodal. A nice exercise for the reader is to check that this condition holds if and only if the equation $p_1 - p_2 + p_3 - p_4 = 0$ holds in the complex plane.

The core $CS_0$ is diffeomorphic to $S^1 \times R^1$.

## 14. Reduced configuration space.

To aid in visualization, and also to help with our proof of the required transversality results, we define the **reduced configuration space** RCS ⊂ CS to be the subset where $p_1 = (1, 0) = 1 + 0i = 1$. Then RCS is homeomorphic to $R^3$ and we use the group structure on $S^1$ to express the homeomorphism $S^1 \times RCS \to CS$ by

$$(e^{i\theta}, (1, p, q, r)) \to (e^{i\theta}, e^{i\theta}p, e^{i\theta}q, e^{i\theta}r).$$

For the purpose of drawing pictures, we change coordinates by writing

$$p = e^{2\pi ix}, \quad q = e^{2\pi iy} \text{ and } r = e^{2\pi iz}.$$

Then RCS $\cong \{(x, y, z) : 0 < x < y < z < 1\}$.

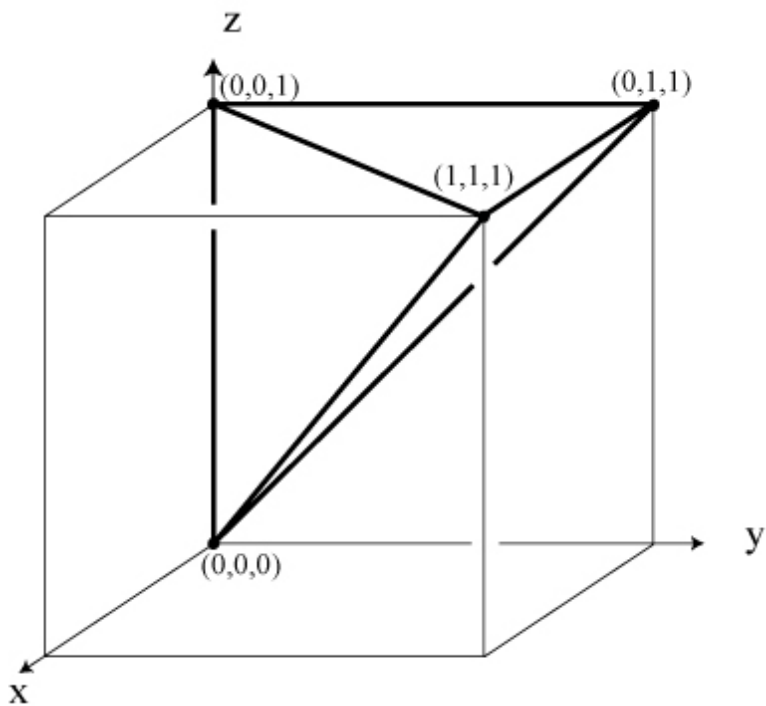

**The reduced configuration space appears as an open solid tetrahedron**

A point $(1, p, q, r)$ in RCS is in the core if and only if $1$ and $q$ are antipodal, and also $p$ and $r$ are antipodal. We let $RCS_0 = RCS \cap CS_0$ denote the core of the reduced configuration space.

In the x, y, z coordinates for RCS, the core $RCS_0$ is given by

$$0 < x < y = \tfrac{1}{2} < z = x + \tfrac{1}{2} < 1,$$

and appears as the open line segment connecting $(0, \tfrac{1}{2}, \tfrac{1}{2})$ to $(\tfrac{1}{2}, \tfrac{1}{2}, 1)$.

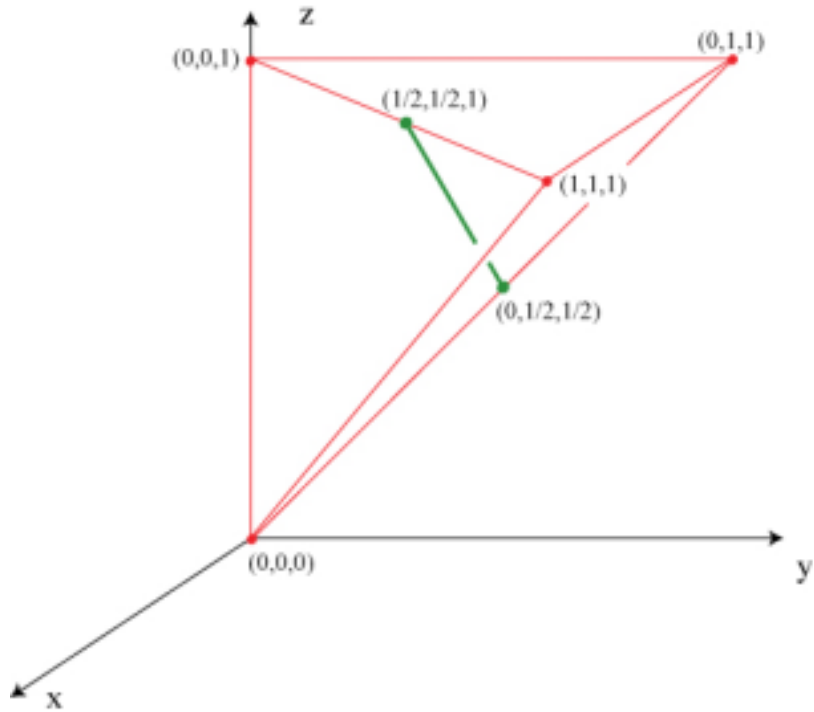

**The core of the reduced configuration space.**

The reduced configuration space has an added advantage which derives from the fact that when we build curves in the plane by cutting the unit circle open at the point $1$, this is already one of the four division points for an element $(1, p, q, r)$ in RCS. This simplifies the construction of the curve, and the computation of its error vector $E(1, p, q, r)$, and we make use of this in the following section.

## 15. The topology of the error map.

The error map $E: CS \to R^2$ takes the core $CS_0$ to the origin, and the complement of the core to the complement of the origin.

Let $\lambda$ be a loop in $CS - CS_0$ which is null-homologous in $CS$ but links the core $CS_0$ once. For convenience, we go down one dimension and show this loop in $RCS - RCS_0$.

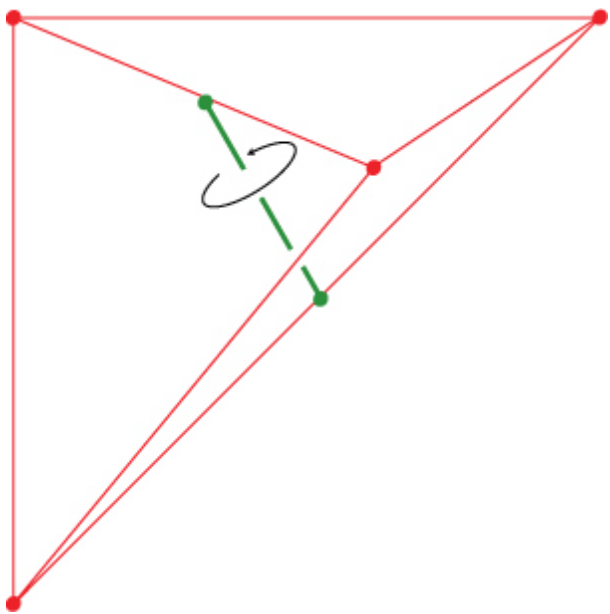

Any two such loops in $CS - CS_0$ are homotopic to one another, up to sign.

**PROPOSITION 15.1.** *The image of the loop $\lambda$ under the error map $E$ has winding number $\pm 1$ about the origin in $R^2$.*

What this proposition says in effect is that plane curves, built with four arcs cut in alternation from two different size circles, are capable of exhibiting the winding number phenomenon which makes the proof of the converse to the Four Vertex Theorem actually work.

Proposition 15.1 is an immediate consequence of the transversality result below.

**PROPOSITION 15.2.** *The differential of the error map* $E: RCS \to R^2$ *is surjective at each point of the core.*

**Proof.** We are given distinct positive numbers $0 < a < b$ and a point $P = (1, p_2, p_3, p_4)$ in RCS. Let $L_1, L_2, L_3, L_4$ denote the lengths of the four arcs $1p_1, p_1p_2, p_2p_3, p_31$ on the unit circle. Scale $a$ and $b$ up or down to new values $a(P)$ and $b(P)$, preserving their ratio, to achieve total curvature $2\pi$:

$$a(P) L_1 + b(P) L_2 + a(P) L_3 + b(P) L_4 = 2\pi.$$

Then construct a curve $\alpha: [0, 2\pi] \to R^2$ according to Dahlberg's plan: begin at the origin and head in the direction of the positive x-axis along circles of curvatures $a(P), b(P), a(P), b(P)$ for lengths $L_1, L_2, L_3, L_4$. The error vector $E(P) = \alpha(2\pi) - \alpha(0)$ indicates the failure to close up.

Consider the figure below, which shows the points $1, p_2, p_3, p_4$ on the unit circle parametrized by arc length $s$, and the corresponding division points $1, q_2, q_3, q_4$ on the unit circle parametrized by angle of inclination $\theta$.

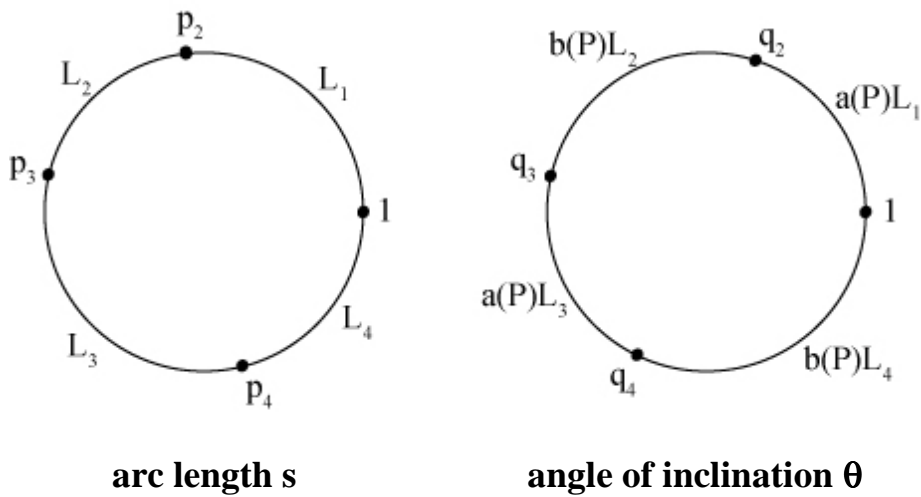

**arc length s**          **angle of inclination $\theta$**

This momentary switch in loyalty from the arc length parameter $s$ to the angle of inclination parameter $\theta$ will be rewarded by an explicit formula for the error map which makes the proof of the current proposition easy.

The map
$$P = (1, p_2, p_3, p_4) \to Q = (1, q_2, q_3, q_4)$$
is a diffeomorphism from the s-parametrized version of RCS to a θ-parametrized version, taking the core $1 - p_2 + p_3 - p_4$ to the core $1 - q_2 + q_3 - q_4$.

This alternative version of RCS helps us to calculate the error map. First note that
$$E(P) = \int_0^{2\pi} e^{i\theta(s)}\, ds = \int_0^{2\pi} e^{i\theta} (ds/d\theta)\, d\theta = \int_0^{2\pi} e^{i\theta} / \kappa(\theta)\, d\theta.$$

This can be computed explicitly as a sum of four integrals, and we get
$$E(P) = \bigl(1/i\, b(P) - 1/i\, a(P)\bigr)\, (1 - q_2 + q_3 - q_4).$$

The differential of this map is particularly easy to compute along the core, since the vanishing of the expression $1 - q_2 + q_3 - q_4 = 0$ there frees us from the need to consider the rate of change of the factor $\bigl(1/i\, b(P) - 1/i\, a(P)\bigr)$.

If we move $q_2$ counterclockwise at unit speed along the circle, we have $dq_2/dt = i\, q_2$, and hence
$$dE/dq_2 = \bigl(1/i\, b(P) - 1/i\, a(P)\bigr)(-i\, q_2).$$

Likewise,
$$dE/dq_3 = \bigl(1/i\, b(P) - 1/i\, a(P)\bigr)(i\, q_3).$$

These two rates of change are independent because $q_2$ and $q_3$ are always distinct, and along the core they can not be antipodal. It follows that the differential of the error map is surjective at each point of the core of the θ-parametrized version of RCS. Since the map between the two versions of RCS is a diffeomorphism taking core to core, the same result holds for our original s-parametrized version of RCS, completing the proof of the proposition.

# 16. Dahlberg's disk **D**.

Dahlberg's choice of 2-cell **D** ⊂ Diff($S^1$) consists of the *special Möbius transformations*

$$g_\beta(z) = (z - \beta) / (1 - \bar\beta z),$$

where $|\beta| < 1$ and $\bar\beta$ is the complex conjugate of $\beta$.

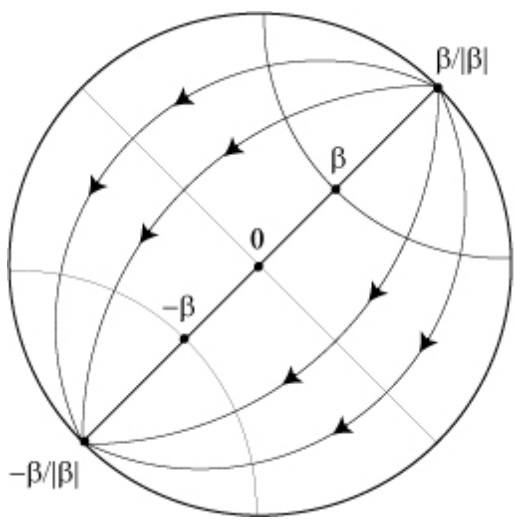

**The action of $g_\beta$ on the unit disk in the complex plane**

These special Möbius transformations are all isometries of the Poincaré disk model of the hyperbolic plane. The transformation $g_0$ is the identity, and if $\beta \neq 0$, then $g_\beta$ is a hyperbolic translation of the line through 0 and $\beta$ which takes $\beta$ to 0 and 0 to $-\beta$. The point $\beta/|\beta|$ and its antipode $-\beta/|\beta|$ on $S^1$ (the circle at infinity) are the only fixed points of $g_\beta$.

A unit complex number $e^{i\theta}$ acts by multiplication on the Poincaré disk model of the hyperbolic plane, and as such is a hyperbolic isometry which rotates the disk by angle $\theta$ about its center. These rotations interact with Dahlberg's transformations $g_\beta$ by the intertwining formula

$$g_{e^{i\theta}\beta}(e^{i\theta} z) = e^{i\theta} g_\beta(z).$$

## 17. The second transversality result.

Let $g_\beta$ be a point in Dahlberg's disk **D**, and $P = (p_1, p_2, p_3, p_4)$ a point in the configuration space CS. We define

$$g_\beta(P) = (g_\beta(p_1), g_\beta(p_2), g_\beta(p_3), g_\beta(p_4)),$$

so that now **D** acts on CS.

**PROPOSITION 17.1.** *The evaluation map* $CS_0 \times \mathbf{D} \to CS$ *defined by* $(P, g_\beta) \to g_\beta(P)$ *is a diffeomorphism.*

**Proof.** The evaluation map $CS_0 \times \mathbf{D} \to CS$ is smooth. To visualize it, start with a point $P = (p_1, p_2, p_3, p_4)$ in the core $CS_0$ and a transformation $g_\beta$ in **D** which takes $P$ to $Q = (q_1, q_2, q_3, q_4)$. In the Poincaré disk, the geodesic with ends at $p_1$ and $p_3$ is a straight line through the origin, and likewise for the geodesic with ends at $p_2$ and $p_4$. The isometry $g_\beta$ takes them to geodesics with ends at $q_1$ and $q_3$, respectively $q_2$ and $q_4$, which appear to the Euclidean eye as circular arcs meeting the unit circle orthogonally, and intersecting one another at the point $g_\beta(0) = -\beta$.

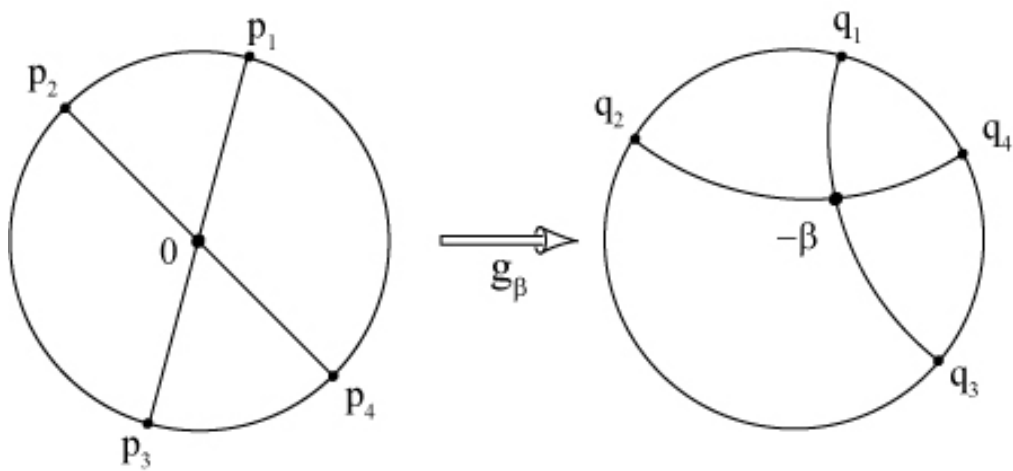

Vice versa, suppose we start with the point $Q = (q_1, q_2, q_3, q_4)$ in CS. Draw the hyperbolic lines, that is, circles orthogonal to the boundary of the disk, with ends at $q_1$ and $q_3$, likewise $q_2$ and $q_4$. The intersection point of these two circular arcs varies smoothly with the four end points.

Call this intersection point $-\beta$. The transformation $g_{-\beta}$ takes the point $-\beta$ to the origin, and hence takes the two hyperbolic lines through $-\beta$ to lines through the origin. It follows that $P = g_{-\beta}(Q)$ lies in the core $CS_0$, and that $g_\beta(P) = Q$. The association $Q \to (P, g_\beta)$ provides a smooth inverse to the evaluation map $CS_0 \times \mathbf{D} \to CS$, showing it to be a diffeomorphism.

**COROLLARY 17.2.** *For each fixed point* P *in the core* $CS_0$, *the evaluation map* $g_\beta \to g_\beta(P)$ *is a smooth embedding of Dahlberg's disk* **D** *into* CS *which meets the core transversally at the point* P *and nowhere else.*

This follows directly from Proposition 17.1.

## 18. The image of Dahlberg's disk in the reduced configuration space.

We seek a concrete picture of Dahlberg's disk **D** in action. Recall from the preceding section that for any point $P = (p_1, p_2, p_3, p_4)$ in CS, the correspondence

$$g_\beta \to (g_\beta(p_1), g_\beta(p_2), g_\beta(p_3), g_\beta(p_4))$$

maps **D** into CS. If we start with the point $P_0 = (1, i, -1, -i)$, and then follow the above map by the projection of CS to RCS, we get the correspondence

$$g_\beta \to \left(1, g_\beta(1)^{-1}g_\beta(i), g_\beta(1)^{-1}g_\beta(-1), g_\beta(1)^{-1}g_\beta(-i)\right).$$

If we convert to xyz-coordinates so that RCS is the open solid tetrahedron introduced in section 14, then the image of **D** is shown below. In the picture on the left, the tetrahedron has its usual position; on the right, it is rotated so as to clearly display the negative curvature of the image disk.

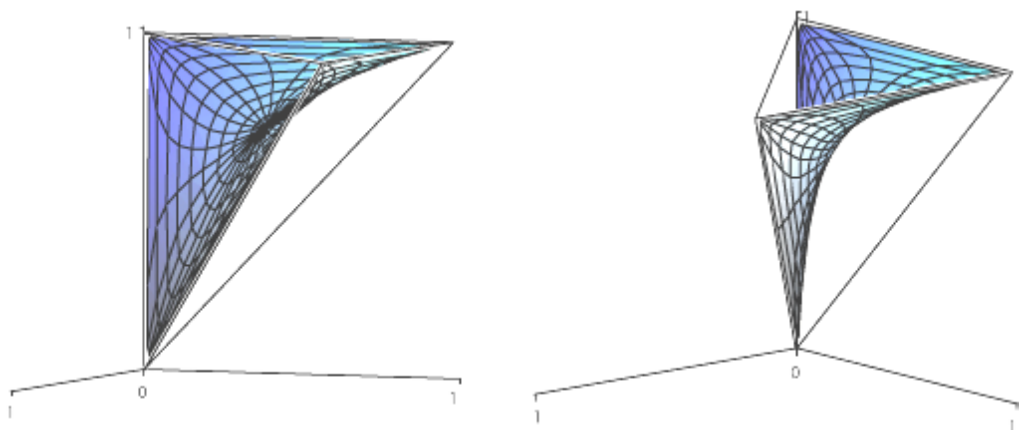

**Dahlberg's disk D mapped into the tetrahedron - two views**

These computer drawn pictures of the image of **D** clearly show that its boundary lies along four of the six edges of the tetrahedron. This is a reflection of the fact that as $\beta$ approaches the boundary of the unit disk, the special Möbius transformation $g_\beta$ leaves the point $\beta/|\beta|$ on the boundary circle fixed, but moves all other points on the boundary circle towards the antipodal point $-\beta/|\beta|$. Hence the four points $1, g_\beta(1)^{-1}g_\beta(i), g_\beta(1)^{-1}g_\beta(-1), g_\beta(1)^{-1}g_\beta(-i)$ limit to at most two points. This "double degeneracy" corresponds to edges (rather than faces) of the tetrahedron.

# 19. Dahlberg's proof.

We are given a continuous curvature function $\kappa \colon S^1 \to R$ which has at least two local maxima and two local minima, and must find an embedding $\alpha \colon S^1 \to R^2$ whose curvature at the point $\alpha(t)$ is $\kappa(t)$ for all $t \in S^1$.

### *Step 1. We temporarily replace $\kappa$ by a curvature step function $\kappa_0$, to which we apply the winding number argument.*

Changing the sign of $\kappa$ if necessary, there are constants $0 < a < b$ so that $\kappa$ takes the values $a, b, a, b$ at four points in succession along $S^1$. Let $\kappa_0$ denote the curvature step function which takes these same values along the four quarter-circles with end points at $1, i, -1, -i$.

Given any $\varepsilon > 0$, we can find a preliminary diffeomorphism $h_1$ of the circle so that the curvature function $\kappa \circ h_1$ is $\varepsilon$-close in measure to $\kappa_0$. We then rescale both $\kappa_0$ and $\kappa \circ h_1$ to have total curvature $2\pi$. They will again be $\varepsilon$-close in measure, for some new small $\varepsilon$.

We apply the winding number argument to the curvature step function $\kappa_0$ as follows.

Consider the point $P_0 = (1, i, -1, -i)$ in the core $CS_0$ of the configuration space $CS$, and map Dahlberg's disk **D** into $CS$ by sending $g_\beta \to g_\beta(P_0)$. By Corollary 17.2, this evaluation map is a smooth embedding of **D** into $CS$ which meets the core $CS_0$ transversally at the point $P_0$ and nowhere else. By Proposition 15.1, each loop $|\beta| = $ constant in **D** is sent by the composition

$$\mathbf{D} \text{ —evaluation map} \to CS \text{ —error map} \to R^2$$

into a loop in $R^2 - $ origin with winding number $\pm 1$ about the origin.

We translate this conclusion into more concrete terms as follows.

Let $c(\beta) \kappa_0 \circ g_\beta$ be the rescaling of the curvature step function $\kappa_0 \circ g_\beta$ which has total curvature $2\pi$, and let $\alpha(\beta) \colon [0, 2\pi] \to R^2$ be the corresponding arc-length parametrized curve with this curvature function. As $\beta$ circles once around the origin, the corresponding loop of error vectors $E(\alpha(\beta))$ has winding number $\pm 1$ about the origin.

*Step 2. We transfer this winding number argument to the curvature function* $\kappa$.

Let $c(h_1, \beta) \, \kappa \circ h_1 \circ g_\beta$ be the rescaling of the curvature function $\kappa \circ h_1 \circ g_\beta$ which has total curvature $2\pi$, and let $\alpha(h_1, \beta) : [0, 2\pi] \to R^2$ be the corresponding arc-length parametrized curve with this curvature function.

Fixing $|\beta|$, we can choose $\varepsilon$ sufficiently small so that each curve $\alpha(h_1, \beta)$ is $C^1$-close to the curve $\alpha(\beta)$ constructed in Step 1. Then as $\beta$ circles once around the origin, the corresponding loop of error vectors $E(\alpha(h_1, \beta))$ will also have winding number $\pm 1$ about the origin. Note that we only need $C^0$-close for this step.

It follows that there is a diffeomorphism $g_{\beta'}$ with $|\beta'| \le |\beta|$ so that $E(\alpha(h_1, \beta')) = 0$, which tells us that the curve $\alpha(h_1, \beta')$ closes up smoothly. If $|\beta|$ and $\varepsilon$ are sufficiently small, then the closed curve $\alpha(h_1, \beta')$ will be as $C^1$-close as we like to the fixed bicircle with curvature $c_0 \kappa_0$, and will hence be simple.

The simple closed curve $\alpha(h_1, \beta')$ realizes the curvature function $c(h_1, \beta') \, \kappa \circ h_1 \circ g_{\beta'}$. Rescaling it realizes the curvature function $\kappa \circ h_1 \circ g_{\beta'}$, and then reparametrizing it realizes the curvature function $\kappa$.

This completes Dahlberg's proof of the Converse to the Four Vertex Theorem.

## 20. Extensions and generalizations of the Four Vertex Theorem.

The Four Vertex Theorem for planar curves may be generalized in the following ways:

   (1) by considering more general surfaces than just the Euclidean plane;
   (2) by relaxing the condition that the curve be simple;
   (3) by considering space curves;  and
   (4) by seeking conditions which would ensure the existence of more than four vertices.

Some generalizations go in several of these directions at the same time.

(1) The consideration of more general surfaces was initiated by Jackson (1945), who showed that any simple closed curve which lies on a surface of constant curvature in Euclidean space, and is null-homotopic on that surface, has at least four vertices, meaning local maxima or local minima of the geodesic curvature. See also Thorbergsson (1976).

(2) The generalization to nonsimple curves was initiated by Pinkall (1987), who showed that if a closed planar curve bounds an immersed surface in the plane, then it must have at least four vertices. He also pointed out that there exist closed curves of any rotation index which have only two vertices.

The generalizations of Jackson and Pinkall may be combined to show that any closed curve on a complete orientable 2-manifold of constant curvature, which bounds an immersed disk, has at least four vertices. See Maeda (1998) and Costa and Firer (2000).

(3) The generalization to space curves stems from the observation that stereographic projection maps the vertices of a planar curve to vertices (with respect to geodesic curvature) of its image on the sphere. This is a special case of the more general phenomenon that any conformal transformation of a surface preserves the vertices. On a sphere, however, the points where the geodesic curvature vanishes correspond to points of vanishing torsion. Thus the fact that a spherical curve has at least four points of vanishing torsion may be considered a generalization of the Four Vertex Theorem in the plane. Sedykh (1994) showed that this result holds for all space curves which lie on the boundary of their convex hull, thus proving a conjecture of Scherk.

(4) The first result seeking more than four vertices is that of Osserman (1985) discussed earlier in this paper. Another generalization was initiated by Pinkall, who had conjectured that if a planar curve bounds an immersed surface of genus $g$, then it must have $4g + 2$ vertices. This was disproved by Cairns, Ozdemir and Tjaden (1992), where the authors constructed for any genus $g \geq 1$ a planar curve bounding an immersed surface of genus $g$ which has only $6$ vertices. Later, Umehara (1994) showed that this result is optimal.

The comments in this section were taken, with almost no changes, from a letter written by Mohammad Ghomi.

# 21. History.

*Syamadas Mukhopadhyaya* was born on June 22, 1866 at Haripal in the district of Hooghly. After graduating from Hooghly College, he took his M.A. degree from the Presidency College, Calcutta, and was later awarded the first Doctor of Philosophy in Mathematics from the Calcutta University, with a thesis titled, "Parametric Coefficients in the Differential Geometry of Curves in an N-space", which earned for its author the Griffith Memorial Prize for 1910.

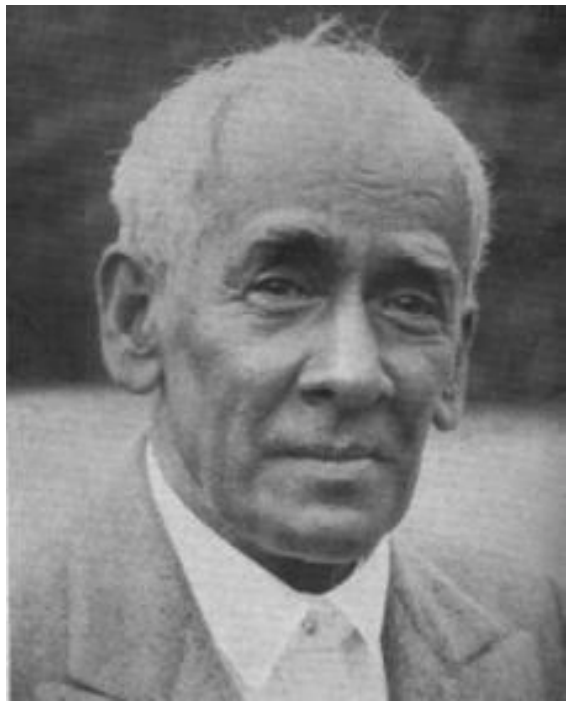

**Syamadas Mukhopadhyaya (1866 - 1937)**

Mukhopadhyaya began work as a Professor in the Bangabasi College, and some years later moved to the Bethune College, both in Calcutta. There he lectured not only on Mathematics, but also on Philosophy and English Literature. He was later transferred to the Presidency College as Professor of Mathematics. When the Calcutta Mathematical Society was started in 1909, he was a founding member, and in 1917 became Vice-President of the Society.

All of Mukhopadhyaya's mathematical works had a strong geometric flavor. He wrote many papers on the geometry of curves in the plane, and his proof of the Four Vertex Theorem in the case of strictly positive curvature came at the very beginning of his research career. His work on the geometry of curves in n-space, the early parts of which were submitted for his doctoral thesis, continued afterwards, and was published in six parts in the Bulletin of the Calcutta Mathematical Society over a period of years.

After his retirement in 1932, Mukhopadhyaya went to Europe to study methods of education, and on returning to India wrote a series of memoirs about this. He was elected President of the Calcutta Mathematical Society, and served in this capacity until his death from heart failure on May 8, 1937.

*"His death has removed from the world of science a man of outstanding genius; but the intimate circle of his friends, admirers and pupils will remember him not only as an able Professor and a successful researcher with a deep insight into the fundamental principals of synthetic geometry but also as a sincere guide, philosopher, and friend, ever ready to help the poor in distress."*

Mukhopadhyaya's obituary appeared in the Bulletin of the Calcutta Mathematical Society, Vol. 29 (1937), pages 115-120. The comments above were taken from this obituary, and the final sentence was copied verbatim.

***Adolf Kneser*** was born on March 19, 1862 in Grüssow, Germany.  He received his Doctor of Philosophy from the University of Berlin in 1884, with a thesis *"Irreduktibilität und Monodromiegruppe algebraischer Gleichungen"* written under the direction of Leopold Kronecker and Ernst Kummer, and influenced by Weierstrass.  Kneser was appointed to the chair of mathematics in Dorpat, and then later to the University of Breslau, where he spent the rest of his career.

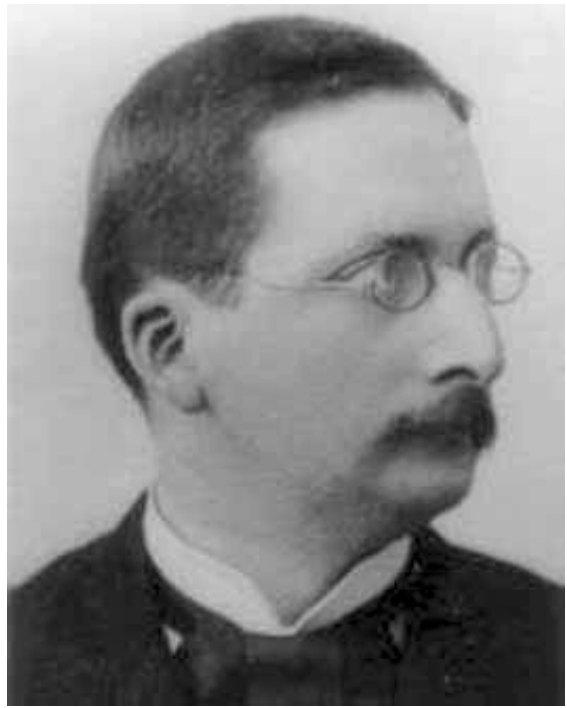

**Adolf Kneser (1862 – 1930)**

After his initial interest in algebraic functions and equations, Kneser turned his attention to the geometry of space curves, and in 1912 proved the Four Vertex Theorem in the general case of mixed positive and negative curvatures.

Adolf's son Hellmuth, born in 1898 while the family was still in Dorpat, followed in his father's mathematical footsteps, and obtained his doctorate from Gottingen in 1921 under the direction of David Hilbert, with a thesis titled, *"Untersuchungen zur Quantentheorie"*.  The next year, in 1922, Hellmuth published a new proof of the Four Vertex Theorem.  Hellmuth helped Wilhelm Süss to found the Mathematical Research Institute at Oberwolfach in 1944, and provided crucial support over the years to help maintain this world famous institution.

Hellmuth's son Martin continued in the family business, receiving his Ph.D. from the Humboldt University in Berlin in 1950 under the direction of Erhard Schmidt, with a thesis titled, *"Über den Rand von Parallelkorpern"*.  He spent the early part of his mathematical career in Munich, and then moved to Gottingen.

Adolf died on January 24, 1930 in Breslau, Germany (now Wroclaw, Poland).
Hellmuth died on August 23 1973 in Tubingen, Germany.
Martin died on February 16, 2004.

The comments above were taken from the online mathematical biographies
at the University of St. Andrews in Scotland, *www.history.mcs.st-and.ac.uk* .

***Björn E. J. Dahlberg*** was born on November 3, 1949 in Sunne, Sweden. He received his Ph.D. in mathematics from Göteborg University in 1971 with a thesis, *"Growth properties of subharmonic functions",* and with Tord Ganelius as advisor. Dahlberg was 21 at the time, the youngest mathematics Ph.D. in Sweden.

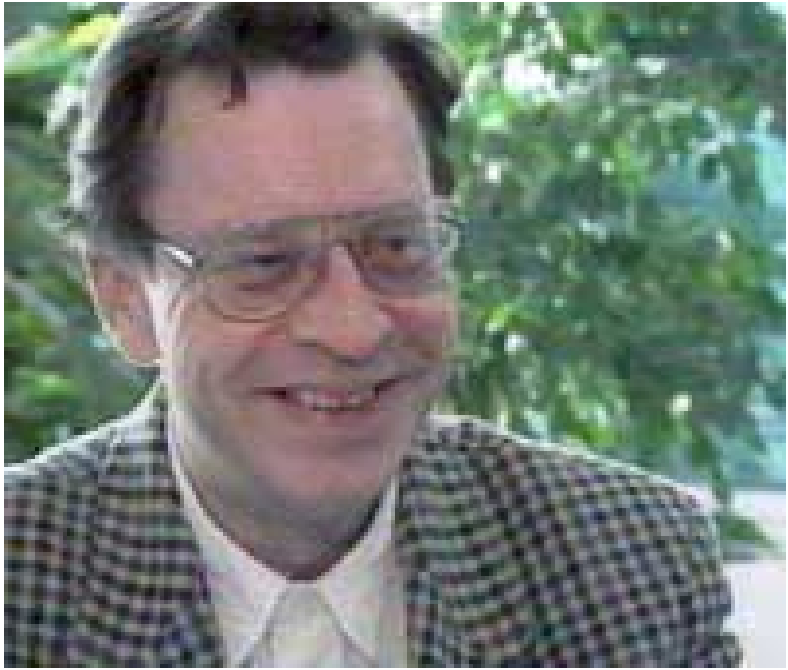

**Björn E. J. Dahlberg (1949 – 1998)**

Dahlberg's mathematical interests were shaped during his postdoctoral fellowship year at the Mittag-Leffler Institute in Stockholm, which Lennart Carleson had recently revitalized. Dahlberg was deeply influenced by Carleson, and participated in the programs of harmonic analysis and quasi-conformal mappings at the Institute, while learning about partial differential equations.

One issue of great importance to the mathematical community at that time was to extend known results in these areas to domains with complicated boundaries. In 1977 Dahlberg showed that, for bounded Lipschitz domains, surface measure on the boundary and harmonic measure are mutually absolutely continuous, which implies the solvability of the Dirichlet problem for the Laplacian on such domains.

For this work, Dahlberg was awarded the Salem Prize (Paris) in harmonic analysis in 1978, the Edlund's prize from the Royal Academy of Science in Stockholm in 1979, and was invited to give a 45 minute address at the International Congress of Mathematicians in Warsaw, Poland in 1982.

In 1980 Dahlberg was appointed to a full professorship at Uppsala University, and three years later to a full professorship at Göteborg, which remained his home base for the rest of his life. He traveled a lot, and held visiting faculty positions at the Universities of Minnesota, Michigan and Texas, at Purdue University and Washington University in St. Louis, at the Universite Paris-Sud, as well as at ETH (Zurich), Yale, Cal Tech, Chicago and South Carolina.

While in Göteborg, Dahlberg started a project with Volvo on computer aided geometric design of the surfaces of a car, leading to new theoretical results and a new design system. In his last years, Dahlberg worked on questions of discrete geometry connected to his work with Volvo and other companies, as well as to his older interest in non-smooth domains.

Dahlberg died suddenly on January 30, 1998 from meningitis, an inflammation of the lining that protects the brain and spinal cord.

The comments above were taken mainly from a letter written by Vilhelm Adolfsson, one of Dahlberg's Ph.D. students who, along with Peter Kumlin, found Dahlberg's manuscript on the Converse to the Four Vertex Theorem after his death, edited it, and submitted it for publication. Ville went on to add some personal notes.

*"Björn was a very inspiring and charismatic person. Where other mathematicians tended to see difficulties, Björn saw possibilities: nothing was impossible! I remember leaving his office after discussing current problems, I always felt encouraged and full of enthusiasm. He often cheerfully exclaimed: The sky is the limit! One might then think that he was overly optimistic, but he was well aware of the difficulties and the uncertainties. One of his favorite sayings was 'Hindsight is the only certain science.' "*

## 22. Acknowledgments.

The entire mathematical community is indebted to Vilhelm Adolfsson and Peter Kumlin for rescuing Dahlberg's proof of the Full Converse to the Four Vertex Theorem, reading and editing it, and shepherding it through to publication. In addition, we are grateful to Ville for sending us a copy of the proof, to both Ville and Peter for telling us about Dahlberg's life and work, and to Ville for the detailed letter which supplied the information for the previous section.

We are enormously grateful to Benjamin Schak and Clayton Shonkwiler, two graduate students in mathematics at the University of Pennsylvania, for reading the entire manuscript and making many improvements, as well as for studying and explaining to us the proofs of Mukhopadhyaya and Kneser.

Special thanks to Mohammad Ghomi for his detailed letter, reproduced with only minor changes in section 20, about extensions and generalizations of the Four Vertex Theorem.

University of Pennsylvania

*deturck@math.upenn.edu*
*gluck@math.upenn.edu*
*dpomerle@math.upenn.edu*
*dvick@math.upenn.edu*